\documentclass[11pt]{article}
\usepackage{amsmath, amssymb, amsfonts, amsthm}

\theoremstyle{plain}
\newtheorem*{theorem}{Theorem}

  \newcommand{\sig}{\sigma}

  \newcommand{\Hom}{\operatorname{Hom}}  
  \newcommand{\simto}{\overset{\sim}{\rightarrow}}

\newcommand{\Cbb}{\mathbb{C}}
\newcommand{\Pbb}{\mathbb{P}}

 \newcommand{\Ocal}{\mathcal{O}}

\begin{document}
  \title{ Equivariant Analogue of Grothendieck's Theorem for Vector Bundles
on $\Pbb^1$}
  \author{Shrawan Kumar}

  \maketitle
{\it \qquad\qquad To Professor C.S. Seshadri on his Seventieth birth day}
\vskip6ex

The aim of this note is to prove the  following equivariant analogue of
 Grothendieck's theorem, which was asked by 
W. Fulton.

Let $T$ be a complex (connected) torus  acting algebraically on 
$X=\Pbb^1(\Cbb )$.

  \begin{theorem}  Let $E$ be a $T$-equivariant algebraic vector bundle on $X$.  Then
there exist $T$-equivariant line subbundles $L_1,\cdots ,L_m\subset E$ such
that we have a $T$-equivariant vector bundle isomorphism:
   \begin{align*}
L_1\oplus \cdots &\oplus L_m \simto E \\
 \; \searrow &\phantom{\Pbb^1} \qquad\swarrow\\
 \;\;\; \phantom{\searrow} &\Pbb^1\,, 
  \end{align*}
under 
$v_1\oplus \cdots \oplus v_m \mapsto v_1 + \cdots + v_m.$
  \end{theorem}

  \begin{proof}  As is well known (cf. [2; Exercise 6.6, Chapter I]), 
the group of algebraic automorphisms of $X$ can be identified with the 
projective linear group PGL$(2)$ via $\begin{pmatrix} a & b \\ c & d 
\end{pmatrix} \cdot 
[z_1,z_2]=[az_1+bz_2, cz_1+dz_2].$ Let $D\subset $ PGL$(2)$ be the 
standard diagonal subgroup. Since $D$ fixes the point at infinity 
$[1,0]$, the action of $D$ on $X$ lifts to give a $D$-equivariant line 
bundle structure on $\Ocal_X(1)$ (and hence on any $\Ocal_X(n), n\in 
\Bbb Z$). Thus $\Ocal_X(n)$ acquires a $T$-equivariant line bundle 
structure.

By the 
virtue of  Grothendieck's theorem (cf. [1; Th\'eor\`eme 2.1]), decompose 
$E \cong \Ocal_X(n_1) \oplus \cdots \oplus \Ocal_X(n_m)$ with $n_1\geq
n_2 \geq \cdots \geq n_m$.  Tensoring $E$ with $\Ocal_X(-n_1)$ and putting a
$T$-equivariant line bundle structure on $\Ocal_X(-n_1)$, we can assume
that $E\simeq \Ocal_X\oplus\Ocal_X(n_2)\oplus \cdots \oplus\Ocal_X(n_m)$
with each $n_2, \cdots ,n_m \leq 0$.
Thus any $0\neq \sig\in H^0(X,E)$ is nowhere vanishing.  Further,
$H^0(X,E)\neq 0$.  Of course, $H^0(X,E)$ is a $T$-module.  Take a
$T$-eigenvector $0\neq \sig\in H^0(X,E)$.  Then $\{ \Cbb\sig
(x)\}_{x\in\Pbb^1}$ is a $T$-equivariant line subbundle, denoted $L_1$, 
 of $E$. 

We get an exact  sequence of $T$-equivariant bundles:
  \begin{equation}
0\to L_1 \to E \to E/L_1 \to 0 .  \tag{$*$}
  \end{equation}
We claim that this sequence is $T$-equivariantly split:

Consider the $T$-module  $\Hom_{\Ocal}(E/L_1,E)$ of bundle morphisms.  Then we
get a natural $T$-module  map
  \[
\pi: \Hom_{\Ocal}(E/L_1,E) \twoheadrightarrow \Hom_{\Ocal}(E/L_1,E/L_1).
  \]
Observe that all the line bundles, occurring in $E/L_1$ as a direct
summand, have degrees $\leq 0$. This can be seen, e.g., by tensoring the 
sequence ($*$) with $\Ocal_X(-1)$ and then considering the associated long 
exact cohomology sequence. Thus ($*$) splits  as vector bundles
by [2; Propositions 6.3 and 6.7, Chap. III]
(without regarding the $T$-equivariance).  In particular, $\pi$ is
surjective.  Thus $\pi$ induces a surjective map
  \[
\pi^T:  \Hom_{\Ocal}(E/L_1,E)^T \twoheadrightarrow
\Hom_{\Ocal}(E/L_1,E/L_1)^T,
  \]
where the superscript $T$ means $T$-invariants.  Take a preimage $f$ of
  the identity homomorphism
$\mathbf{I}$ under $\pi^T$.  This $f$ provides a $T$-equivariant
splitting of ($*$).  Thus $E \simeq L_1\oplus E/L_1$ as $T$-equivariant
bundles.  So, by induction on rank $E$, the theorem follows.
   \end{proof}

\vskip3ex
\noindent
{\bf Acknowledgements.}  
The author was supported by NSF.  
\vskip6ex
\noindent
{\bf References.}

\noindent
[1] Grothendieck, A.: Sur la classification des fibres holomorphes sur la 
sphere de
Riemann, {\it Amer. J. Math. 79} (1957), 121--138.

\noindent
[2] Hartshorne, R.: {\it Algebraic Geometry}, GTM 52, Springer-Verlag, 
1977.
\vskip4ex
\noindent
{\bf Address:}

\noindent
Department of Mathematics, UNC at Chapel Hill, Chapel Hill,
NC 27599--3250, USA
 \end{document}